\RequirePackage{fix-cm}

\documentclass[11pt,twoside,a4paper]{article}

\usepackage{amsthm}
\usepackage{amsmath}
\usepackage{amssymb}
\usepackage{graphicx}

\newtheorem{theorem}{Theorem}[section]

\theoremstyle{definition}
\newtheorem{definition}[theorem]{Definition}

\theoremstyle{remark}
\newtheorem{example}[theorem]{Example}

\newcommand{\B}{{\rm Bel}}

\newcommand\restr[2]{{
  \left.\kern-\nulldelimiterspace 
  #1 
  \vphantom{\big|} 
  \right|_{#2}
  }}

\mathchardef\mhyphen="2D
 
\numberwithin{equation}{section}

\begin{document}
 
\title{Assessing forensic evidence by computing belief functions}
\author{Timber Kerkvliet and Ronald Meester \\ \,\, \\ VU University Amsterdam}
\maketitle

\begin{abstract} 
We first discuss certain problems with the classical probabilistic approach for assessing forensic evidence, in particular its inability to distinguish between lack of belief and disbelief, and its inability to model complete ignorance within a given population. We then discuss Shafer belief functions, a generalization of probability distributions, which can deal with both these objections. We use a calculus of belief functions which does not use the much criticized Dempster rule of combination, but only the very natural Dempster-Shafer conditioning. We then apply this calculus to some classical forensic problems like the various island problems and the problem of parental identification. If we impose no prior knowledge apart from assuming that the culprit or parent belongs to a given population (something which is possible in our setting), then our answers differ from the classical ones when uniform or other priors are imposed. We can actually retrieve the classical answers by imposing the relevant priors, so our setup can and should be interpreted as a generalization of the classical methodology, allowing more flexibility. We show how our calculus can be used to develop an analogue of Bayes' rule, with belief functions instead of classical probabilities. We also discuss consequences of our theory for legal practice.
\end{abstract}


\section{Introduction and motivation}
\label{introduct}

It has been debated for several decades as to what extent theories of probability, are useful and/or suitable for assessing the value of evidence in legal and forensic settings, see e.g.\ \cite{fried96}, \cite{dawid}, \cite{rob_vig95}, \cite{ait95}, \cite{sjerps00}. The debate mainly concentrates on the question whether or not the classical theory of probability, by which we mean the theory following Kolmogorov's axioms, is suitable in legal problems. 

The current dominant view proposes that we should use the classical probability axioms in court, in particular the axiom of additivity, i.e.\ $P(A) + P(B) = P(A \cup B)$ whenever $A$ and $B$ are disjoint (\cite{fried96}, \cite{dawid}, \cite{rob_vig95}, \cite{ait95}, \cite{sjerps00}).  We are typically interested in an event, often denoted by $G$, that a given individual is the donor of a DNA profile, is the criminal in a certain crime, is the father of a certain child, or likewise. The main tool in this dominant view is Bayes' formula
\begin{equation}
\label{bayesintro}
\frac{P(G|E)}{P(\bar{G}|E)}= \frac{P(E|G)}{P(E|\bar{G})} \cdot \frac{P(G)}{P(\bar{G})},
\end{equation}
where $E$ denotes relevant evidence, and where $\bar{G}$ denotes the complement of $G$ (or another set such that $G \cap \bar{G}=\emptyset$). 

Bayes' formula transforms prior odds $\frac{P(G)}{P(\bar{G})}$ into posterior ones $\frac{P(G|E)}{P(\bar{G}|E)}$ by multiplying the prior odds by the likelihood ratio $\frac{P(E|G)}{P(E|\bar{G})}$. We will give more detail on the legal practice and the use of Bayes' rule in Section \ref{practice}, but note that the use of Bayes' formula presupposes the idea that all quantities of interest, including the prior, can indeed be expressed as probabilities satisfying the usual Kolmogorov axioms. 

The alternative view insists that classical probability theory, with the axioms of Kolmogorov, is in many cases not suitable to be used in court or in forensics, for various reasons. We adhere to this alternative view, since we believe that there are situations in which the axioms of Kolmogorov are too restrictive, and we start by giving a number reasons and examples to support this claim. 

First, it has been observed by many that the classical theory cannot distinguish between lack of belief and disbelief. Here, disbelief is associated with evidence indicating the negation of a proposition, whereas lack of belief is associated with not having evidence at all. As Shafer \cite{shafer76} puts it, the classical theory does not allow one to withhold belief from a proposition without according that belief to the negation of the proposition. When we want to apply a theory of probabilities to legal issues, this becomes a relevant issue. Indeed, if certain exculpatory evidence in a case is dismissed, then this may result in less belief in the innocence of the suspect, but it gives no further indication for guilt. 

The second shortcoming of the classical theory is its inability to model complete ignorance within a given population. We first give two examples, and then elaborate on this issue.

\begin{example}
\label{ex1a}
In The Netherlands, a well known court case concerned a traffic accident caused by a car with two passengers. Although it was not disputed that the car caused the accident, it was unclear which of the two passengers was driving. The classical solution to deal with this, is to impose a fifty-fifty prior on the two passengers, but this is in fact not corresponding to reality. In reality we know that one of the two passengers drove, but we are otherwise ignorant. This cannot be modeled with classical probability.     
\end{example}

This example can  be generalized into the well known and classical island problem:

\begin{example}
\label{ex2}
In the classical version of the island problem (see e.g.\ \cite{sloomee11} and \cite{baldon95}) a crime has been committed on an island, making it a certainty that an inhabitant of the island committed it. In the absence of any further information, the classical point of view is to assign a uniform prior probability over all inhabitants concerning the question who is the culprit. The combination of assigning probability 1 to the collection of all inhabitants and probability 0 to each individual is impossible under the classical axioms of probability, although this may be exactly the prior one needs and wants to impose.
\end{example}

The last two examples may need some elaboration since it may not be so obvious why ignorance cannot be properly modeled by a uniform distribution over all possibilities. 

Firstly, when we look at the island problem in Example \ref{ex2}, it is simply the case that we do not have information pointing to any individual. We do have group information, but no individual information.  With a uniform distribution over the group, you nevertheless make a statement about each individual. This is very relevant in legal cases since these are against individuals, not against a whole population. 

Secondly, it is simply not the case that a uniform distribution does not convey any information. Even in a frequentistic context, a uniform distribution tells us something when we repeat the experiment many times. Or, to phrase the same point differently, having probability $\frac12$ for head to come up in a coin flip, is information.  

Finally, an uninformative prior leads to different results than a uniform prior in our theory, as we will see in Section \ref{island}. The fact that these priors lead to different results confirms that these priors are really distinct: a uniform prior is not a prior representing ignorance, and using a uniform prior does not lead to the same results as using a prior that does represent ignorance. 

The examples suggest that the usual axioms of probability may not always be appropriate in legal and forensic settings. In particular, there are at least two (related) problems: (1) the additivity of probabilities, that is, $P(A) + P(B)= P(A \cup B)$ if $A$ and $B$ are disjoint, is not always desirable, and (2) there is no way to model ignorance in a given population in the classical theory. We are not the first ones to observe this, of course. Already back in the seventies of the previous century, there have been at least two major attempts, by Cohen \cite{cohen77} and Shafer \cite{shafer76} respectively, to develop a theory of probabilities, or generalizations thereof, for legal settings outside the realm of the axioms of Kolmogorov. However, both these attempts have been criticized fiercely (references below), for various and different reasons, and nowadays they are not used at all in assessing evidence in legal settings. 

This history fits in a classical pattern in science at large. There is a certain theory (in our case classical probability axioms) which is supposed to describe or explain certain phenomena (in our case dealing with uncertainty in legal and forensic context). Then certain problems or anomalies arise (for our case, see the examples above). Nevertheless, despite these anomalies, the theory is often upheld, often mainly by the lack of an acceptable alternative. In our case, the proposed alternatives of Cohen and Shafer seemed to have too many problems in themselves, and as such they were not acceptable and the classical theory prevailed. The philosopher Thomas Kuhn described this generic process for instance in \cite{kuhn96}, together with many examples from the history of science.

We think that Shafer's approach is the most promising when it comes down to application possibilities and general acceptance, since it is conceptually simpler and closer to classical probability than the approach of Cohen. This is the reason that in this article we restrict our attention to Shafer's approach. 

Before we introduce belief functions properly, we should review for what reasons Shafer's belief functions have been essentially ignored in the forensic and mathematical literature, other than being criticized.  To be sure, they are not completely ignored, but they are typically framed as somewhat of a curiosity and not taken very seriously, for instance in Dawid's otherwise excellent notes \cite{dawid}. We list three main points of concern that can be found in the literature, and briefly indicate our position.

First of all, Shafer himself reports in \cite{shafer81} that many of his critics rejected his belief functions because of the lack of a suitable betting interpretation. In other words, it is the interpretation that seemed to be an obstacle. However, this criticism does not target the theory, but only the development of the theory. Only if it turns out that a suitable betting intepretation is not possible, it would be a basis to reject the theory. In this light, this criticism should only be seen as a request to further solidify the underpinnings of the theory by a betting interpretation. Shafer argues in \cite{shafer81} that no such behavioral interpretation is necessary. We do think it is a legitimate request and note that in our companion paper to the current one \cite{kerkmee15}, we do formulate a very natural betting interpretation of Shafer's belief functions.

A second, and much more important point of concern about Shafer's belief functions can be found for instance in \cite{schum94} and \cite{fs86}. In both references, the main reason to reject Shafer's belief functions is that the calculus of these belief functions, as put forward in the so called Dempster rule of combination, is arbitrary, not well founded and therefore unacceptable. This rule is supposed to describe how different belief functions should be combined into a new one. In the current article, however, we do not use Dempster's rule of combination. The only thing we need is the Shafer-Dempster conditioning, which we motivate without deriving it from Dempster's rule of combination and is much less controversial, if at all. Hence the current article, is consistent with any stance on Dempster's rule of combination. We do note, however, that we do have an opinion on the matter; we in fact reject Dempster's rule, and motivate this in our companion paper \cite{kerkmee15}. 

A third point of concern is articulated, again, in \cite{fs86}. The very fact that a belief function can allow zero belief to both an event and its complement, makes it, according to \cite{fs86}, inadequate to be used in legal matters, where a decision has to be taken, and where not making a decision is not an option. However, we do not think that this objection is well founded, since it seems to mix up the notion of belief with the act of making a decision about declaring someone guilty. Based on a belief function which assigns belief zero to both guilt and innocence, a suspect will not be convicted. Actually, belief functions seem to do more justice to the situation, since a judge will make a certain decision only if there is enough evidence. One should only convict someone if the belief in the hypothesis that he or she is actually guilty, is high enough, and this does not seem to have anything to do with the fact that the belief in a proposition and its complement can both be zero.

Finally, in \cite{pearl} it is questioned whether or not belief functions respect some `rules' of reasoning when it comes to knowledge. For the most part, this criticism does not apply to our theory and how we want to use it, and we hope to convince the reader that there are many situations in which the use of belief functions is very reasonable, by supplying a number of examples.

The goal of the present article is to convince the forensic and legal communities of the fact that Shafer's belief functions can be put to good use in legal and forensic matters, by using a calculus without Dempster's rule, thereby taking away the main obstacle for the use of Shafer's theory and providing the communities with an acceptable extension of the classical theory. We do this by computing (conditional) belief functions in classical problems like the island problem and parental identification problems. Furthermore, we show through examples how we can use an analogue of Bayes' formula. Shafer's belief functions can take care of the problems noticed in the two examples above. Indeed, they are so general that they can distinguish between lack of belief and disbelief, and they are also flexible enough to be able to model complete prior individual ignorance.

We want to stress that the belief functions of Shafer are a generalization of classical probabilities, and that everything that can be modeled with classical probability theory, can therefore also be modeled with these belief functions. If there is a good reason to take a classical informative prior, the theory allows for that. We lose nothing. As expected, having no prior information leads to different outcomes. This is only reasonable we think, since the classical procedure imposes a prior which is legally problematic.

Although we do address and resolve the problems with the classical theory that we mentioned above, this does not mean that the theory we present solves all problems. The difficulty of quantifying evidential support is just one example of a problem we are still facing in our theory, as is the choice of the relevant population. Despite this, we think the theory is a significant step forward, and perhaps it is hard to imagine how a mathematical theory would satisfactorily solve these issues anyway.

On the theoretical part we will in this article be as brief as possible, but we make sure that the current paper remains self-contained. An in-depth and extensive theoretical development is carried out in the companion paper to the current paper in \cite{kerkmee15}, including an extensive discussion of betting interpretations, a law of large numbers for belief functions, a thorough discussion of independence, a detailed analysis of Dempster's rule, and an in-depth discussion of the interpretation of belief functions.  For the forensic and legal applications in the current paper, such an in-depth mathematical study is not necessary, but we do note the fact that the subject is very interesting from both a theoretical and an applied perspective, and we plan to develop it much more in the future.  

The current paper is organized as follows. In Section \ref{twee} we discuss the basic theory by introducing belief functions and conditioning. Next we apply the theory in Section \ref{island} to the classical island problems, and in Section \ref{parental} to the problem of parental identification. In Section \ref{rule} we explain our analogue of Bayes' rule, and finally in Section \ref{practice} we conclude by discussing some consequences of our theory for legal and forensic casework.  

\section{The basic theory}
\label{twee}

Let $\Omega$ be a finite outcome space, for instance the members of a certain population. We want to make statements about the elements of $\Omega$ in the presence of uncertainty, like saying something about who is the culprit in a certain crime. The classical way to do this 
is by means of a suitable probability distribution on $\Omega$. A probability distribution assigns a non-negative support $p(\omega)$ to each element $\omega \in \Omega$ in such a way that the total support is equal to 1. We may, for instance, express our uncertainty about who is the culprit by means of such a probability distribution. The probability that the culprit can be found in a subset $A$ of $\Omega$ is then equal to 
\begin{equation}
\label{kansen}
P(A):=\sum_{\omega \in A} p(\omega).
\end{equation}
The probability measure $P$ can be interpreted as subjective, frequentistic or otherwise, depending on the context and personal 
taste. The support $p(\omega)$ represents the probability or confidence in the outcome $\omega$, and $P(A)$ represents our probability or confidence in an outcome which is contained in $A$. In classical probability theory, a subset of $\Omega$ is also called an {\em event} or sometimes also a {\em hypothesis}, we make no distinction between the two phrases. The probability measure $P$ describes the probability of all such events or hypotheses.

Next we define basic belief assignments and belief functions. The difference between a {\em basic belief assignment} and a probability distribution, is that the former assigns support to nonempty \emph{subsets} of $\Omega$ rather than to individual outcomes. We write $2^{\Omega}$ for the collection of all subsets of $\Omega$.

\begin{definition}
A function $m: 2^{\Omega} \rightarrow [0,1]$ is a {\em basic belief assignment} if $m(\emptyset)=0$ and
\begin{equation}
\label{eq:additivitym}
\sum_{C \subseteq \Omega} m(C)=1.
\end{equation}
\end{definition}
Whereas $p(\omega)$ represents the probability or confidence in the outcome $\omega$, $m(C)$ represents our confidence in an outcome in $C$ which is not specified further. It may appear that there is not much difference between $P$ and $m$, but in fact there is. The crucial difference between $P$ and $m$ is that the support of a subset $C$ of $\Omega$  is not immediately related to the support of the elements or subsets of $C$. 
For instance, if we have no clue whatsoever about the outcome, that is, if we have no information at all, then we may express this by putting $m(\Omega)=1$ and $m(A)=0$ for all strict subsets of $\Omega$. Or we may take $m(\{a,b\})=1/2$ and simultaneously $m(\{a\})=m(\{b\})=0$, for $a, b \in \Omega$, meaning that we have evidence for the union of $a$ and $b$, but no further information to distinguish between them. 

It is also possible that a basic belief functions only assigns positive support to singletons. In such a case, we are back in the classical situation. The quantity $m(C)$ is sometimes referred to as the {\em evidential support} of $C$.
We should view $m$ as the analogue of $p$ in the classical description above. 

We next define the analogue of $P$, which is called a {\em belief function}. We want to quantify how much belief we can assign to a subset $A$ of 
$\Omega$. To this end, we consider all sets $C$ in $\Omega$ with $C \subseteq  A$, which are precisely the events whose occurrence implies the occurrence of $A$.
The belief in a set $A$ now is the sum of the support of all subsets of $A$.  In terms of evidence, the belief in $A$ is the total evidential support of everything implying $A$. 

\begin{definition}
\label{def:belief}
Given a basic belief assignment $m: 2^{\Omega} \rightarrow [0,1]$, the corresponding {\em belief function} $\B: 2^\Omega \rightarrow [0,1]$ is defined by
\begin{equation}
\B(A) := \sum_{C \subseteq A} m(C).
\end{equation}
\end{definition}

We next discuss a number of examples which should convince the reader that in many situations there are natural belief functions which adequately describe the situation. 

\begin{example}
\label{exsecond}
Suppose we want to state our beliefs about a suspect being guilty or innocent, so $\Omega = \{$guilty, innocent$\}$ is our outcome space. If the evidential support of the the suspect being innocent is $p$, and we have no further information, then we have $m(\{$innocent$\})=p$, $m(\{$guilty$\})=0$ and $m(\Omega)=1-p$. The support of $1-p$ assigned to $\Omega$ should be interpreted as the amount of ignorance. Notice that the belief that the suspect is guilty is not equal to 1 minus the belief that the suspect is innocent. The corresponding belief function $\B$ is given by $\B(\{$guilty$\})=0$, $\B(\{$innocent$\})=p$ and $\B(\Omega)=1$.
\end{example}

\begin{example}
\label{exfirst}
The function $m$ for which $m(\Omega)=1$ and $m(A) =0$ for all other $A \subseteq \Omega$ is a basic belief assignment. The corresponding belief function assigns belief 1 to $\Omega$ and belief zero to all strict subsets of $\Omega$. This belief function expresses total ignorance, except for the fact that the outcome must be in $\Omega$. As such it addresses the problem noticed in Example \ref{ex2}.
\end{example}

\begin{example} ({\em Probability distributions})
Every probability distribution is a belief function, as we already indicated. To see this, let $P: 2^\Omega \rightarrow [0,1]$ be a probability distribution. Set $m(\{ \omega \})=P(\{ \omega\})$ for all $\omega \in \Omega$ and $m(C)=0$ for all $C$ such that $|C|>1$. Then we get
\begin{equation}
\B(A) = \sum_{a \in A} m(\{a\}) = \sum_{a \in A} P(\{a\}) = P(A)
\end{equation}
for every $A \subseteq \Omega$. Probability distributions are belief functions for which the corresponding basic belief assignment only assigns positive support to singletons. If $m(C)>0$ for some $C$ with $|C|>1$, then $\B$ is not a probability distribution because it not additive: for any nonempty, disjoint $A,B \subseteq \Omega$ such that $A\cup B = C$ we find
\begin{equation}
\B(A \cup B) > \B(A) + \B(B).
\end{equation}
\end{example}

In concrete situations, the basic belief assignment can very well be based on classical probabilistic considerations when it is reasonable to do so. Here is an important example which we will discuss in full detail in Section \ref{island}.

\begin{example} ({\em The island problem})
\label{exthird}
Let $X=\{1,\ldots,N+1\}$ be the population of the island. At the scene of the crime a DNA profile $\Gamma$ is found which we know has frequency $p \in (0,1]$ in the population. This means that a randomly chosen person has probability $p$ to have the characteristic, independent of the other individuals. We remark that this assumption is in the realm of classical probability theory. This is reasonable since the frequency interpretation of classical probability works well within the context of DNA profiles. 

Our basic belief assignment should capture our prior knowledge, that is, prior to the fact that we found the DNA profile at the crime scene. We have prior knowledge about two different things: (1) we know that we have selected $S$ uniformly at random from the population, and (2) we know the population frequency of the DNA profile to be $p$. Both these items can be satisfactorily described with a classical probability distribution, and with classical independence assumptions. This leads to the following basic belief assignment, prior to the evidence.  

We set\footnote{We write $A \times B \times C$ for the collection of triples $(a,b,c)$ with $a \in A, b \in B$ and $c \in C$.} $\Omega = X \times X \times \{0,1\}^{N+1}$ and let $C: \Omega \rightarrow X,S: \Omega \rightarrow X$ and $\Gamma_i: \Omega \rightarrow \{0,1\}$ be projections on respectively the first, second and $i+2$-th coordinate. $C$ represents the criminal, $S$ the selected individual, and  $\Gamma_i=1$ indicates that the $i$-th individual has characteristic $\Gamma$. Without any reference to the crime, but with reference to the particular characteristic $\Gamma$, we can model the characteristics $\Gamma_i$ and the choice of $S$ by defining the following basic belief assignment on $\Omega$. Let $y_1, \ldots, y_{N+1} \in \{0,1\}$ be such that $\sum_{i=1}^{N+1} y_i=k$. Then for $k=0, \ldots, N+1$,
\begin{equation}
\label{basis}
m(C \in X, S=x, \Gamma_1 = y_1 , \ldots , \Gamma_{N+1}= y_{N+1}) = \frac{1}{N+1} p^{k} (1-p)^{N+1-k},
\end{equation}
for all $x \in X$, and $m(A)=0$ for all other subsets of $\Omega$. We will typically write $\{S=x, \Gamma_1 = y_1 , \ldots , \Gamma_{N+1}= y_{N+1}\}$ for the set in (\ref{basis}) and similar ones later on, and not explicitly mention the first coordinate. 

This basic belief assignment expresses the facts that characteristics are random and independent with success probability $p$ and that $S$ is chosen uniformly on $X$. But also, and very importantly, it says nothing about the crime and nothing about the identity of the criminal $C$. For instance, for all $x \in X$, we assign zero belief to the event that $C=x$.
\end{example}

In the theory we are about to develop, classical probability distributions are replaced by the more general belief functions, allowing for more flexibility. In Examples \ref{exfirst} and \ref{exthird} above, the basic belief assignments reflect prior knowledge, or the absence thereof. As such it is reasonable to call the basic belief assignments and the corresponding belief functions our {\em priors} in these examples. These priors play the same role as the prior probabilities in the classical situation, with the crucial difference that they need not be a classical probability distribution anymore. As such they need not be additive and they can be genuinely uninformative if that is what corresponds to reality.

\subsection{Conditioning}

The next item on the agenda is to investigate how belief functions change when additional or new information is provided. This is akin to the classical situation in which a prior probability is updated into a posterior one as we briefly mentioned in the introduction. We explain how this works in our setting.

The rule we propose for conditioning is described as follows. Suppose we have a basic belief assignment $m$ and corresponding belief function $\B$. We want to condition on an event $H$. The evidential support $m(A)$ of $A$ now becomes evidential support of $A \cap H$ if $A$ is consistent with $H$ in the sense that $A \cap H \not=\emptyset$. If $A \cap H =\emptyset$, then the new evidential support of $A$ becomes zero. Next we rescale the support in such a way that the support again sums up to $1$. This can of course only be done if the belief in $H^c$ is not $1$.  This leads to the following definition.

\begin{definition}
\label{def:conditionalm}
Let $m: 2^\Omega \rightarrow [0,1]$ be a basic belief assignment and $\B$ the corresponding belief function. For $H \subseteq \Omega$ such that $\B(H^c) \not=1$  we define the {\em conditional basic belief assignment} $m_H: 2^\Omega \rightarrow [0,1]$ by
\begin{equation}
\label{eq:mcond}
m_H(A) := \frac{\sum_{B \cap H = A} m(B)}{1 - \sum_{B \cap H = \emptyset} m(B)},
\end{equation}
for $A \not= \emptyset$ and $m_H(\emptyset)=0$.
The corresponding {\em conditional belief function} $\B_H: 2^\Omega \rightarrow [0,1]$ is defined as
\begin{equation}
\label{conditionalbelief}
\B_H(A)  =  \sum_{B \subseteq A } m_H(B). 
\end{equation}
\end{definition}

The notion of conditioning in (\ref{eq:mcond}) is known as \emph{Dempster-Shafer} conditioning. Shafer \cite{shafer76} derives the formula as a special case of Dempster's rule of combination. We, however, see Definition \ref{def:conditionalm} as standing on its own, motivated by the description above. 

In the special case that $\B=P$ is a probability distribution, our notion of conditional belief coincides with the notion of conditional probability. Indeed, in that case we can write
\begin{equation}
\B_H(A) = \frac{\sum_{\omega \in A \cap H} m(\{\omega\}) }{\sum_{\omega \in H} m(\{\omega\})} = P(A|H),
\end{equation}
for every $A \subseteq \Omega$ and $H$ such that $1-\B(H^c)=P(H)>0$. 

\begin{example}
\label{ex3}
Suppose we have a case in which the suspects are two parents and their son, so $\Omega=\{\mathrm{Father},\mathrm{Mother},\mathrm{Son}\}$. We have a lot of evidence that points to the parents, none of which points to one of them in particular. Further, we have some evidence that points to the son. The corresponding basic belief assignment is, say, $m(\{\mathrm{Father},\mathrm{Mother}\})=\frac{9}{10}$ and $m(\{\mathrm{Son}\})=\frac{1}{10}$. Under the hypothesis $H$ that it is a man, i.e. $H=\{\mathrm{Father},\mathrm{Son}\}$, the evidence pointing to the parents counts as evidence pointing to the father, so
\begin{equation}
\label{condi}
 m_H(\{ \mathrm{Father} \}) = \frac{9}{10}.
\end{equation}
\end{example}

\section{The island problems}
\label{island}

Now that we have established and discussed the main theoretical issues, we move on to our first important examples, namely the well known and classical island problems. The context of the island problems is classical. A crime has been committed on an island with $N+1$ inhabitants, so that we can be sure that one of them is the culprit. Now some characteristic of the criminal, e.g.\ a DNA profile, is found at the scene of the crime and we may assume that this profile originates from the culprit. Then we somehow select an individual $s$ from the island, who happens to have the same characteristic as the criminal. The question is what we can say about the probability or belief in the event that $s$ is in fact the criminal. This is not a well defined question yet, as it depends on the way $s$ was found. We distinguish between two cases: the cold case in which we randomly select an inhabitant, and the search variant in which we consider the inhabitants one by one in a random order, until we find an inhabitant with the characteristic found at the crime scene.

With the island problems, our belief functions allow us to assign a zero prior belief in the guilt of any of the individuals of the island, while at the same time assigning belief one to the full populations. This seems better suited to the legal context than the classical Bayesian setting since assigning a non-zero prior probability to an individual, without any evidence against the individual itself other than belonging to the population, seems unreasonable. Of course, we have to make modeling assumptions, as in the classical case and we will discuss these below. But as we shall see, the outcomes are different from the classical outcome:  if we assume total prior ignorance the belief that $s$ is the culprit is in our setting different from the classical probability that he is guilty under a uniform prior. We turn to the examples now.  

\subsection{The cold case}
\label{cold}

We continue Example \ref{exthird}. We write $A(x, y_1, \ldots, y_{N+1})$ for the set $\{S=x, \Gamma_1 = y_1 , \ldots , \Gamma_{N+1}= y_{N+1}\}$. Now we want to incorporate the crime and condition on the event $H_s=\{S=s , \Gamma_C = \Gamma_S = 1\}$ that $s \in X$ was chosen, and that the criminal $C$ and the chosen $s$ both have characteristic $\Gamma$. According to our theory of conditioning, we have to assign the mass originally assigned to $A(x, y_1, \ldots, y_{N+1})$ in (\ref{basis}) to the intersection of this set with $H_s$, and normalize suitably. Clearly, this intersection is non-empty only if $x=s$ and $y_s=1$, and in this case the intersection $A(s, y_1, \ldots, y_{N+1}) \cap H_s$ can be written as
\begin{equation}
\label{basis2}
\{ C \in \{ i \;:\; y_i=1 \} , S=s; \Gamma_1 = y_1 , \ldots , \Gamma_{N+1}= y_{N+1} \}.
\end{equation}
Next we need to find the correct normalization. We claim that the total mass of sets of the form in (\ref{basis}) which have non-empty intersection with $H_s$ is given by
\begin{equation}
\frac{1}{N+1}\sum_{j=0}^N {N \choose j } p^{j+1}(1-p)^{N-j} = \frac{p}{N+1}.
\end{equation}
Note that one can obtain the outcome at the right hand side by either computing the sum, or by simply noticing that for the intersection to be non-empty it is necessary and sufficient that $s=x$ and that $y_s=1$. Since (\ref{basis}) describes a classical probabilistic experiment, the total mass of the basic belief assignment that has non-empty intersection with $H_s$, is just the probability that in this classical experiment, $S=s$ and $y_s=1$. This happens with probability $p/(N+1)$.

It follows that for $y_1, \ldots, y_{N+1}$ such that $y_s=1$ and $\sum_i y_i = k+1$, $k=0,\ldots, N$ that
\begin{eqnarray*}
m_{H_s}(A(s, y_1, \ldots, y_{N+1})\cap H_s) &=& \frac{p^{k+1}(1-p)^{N-k}}{N+1} \cdot \frac{N+1}{p}\\
&=& p^k(1-p)^{N-k}.
\end{eqnarray*}
Summarizing, we have that for $y_1, \ldots, y_{N+1}$ such that $y_s=1$ and $\sum_i y_i = k+1$,
\begin{equation}
\label{uiteindelijk}
m_{H_s}(C \in \{ i \;:\; y_i=1 \} , S=s; \Gamma_1 = y_1 , \ldots , \Gamma_{N+1}= y_{N+1})=p^k(1-p)^{N-k}.
\end{equation}
Another way of writing this is as follows: for any $A \subset X$ such that $s \in A$ and $|A|=k+1$, we have, writing $1_A(i)$ for $1_{\{i \in A\}}$,
\begin{equation}
m_{H_s}( C \in A , S=s, \Gamma_1 = 1_A(1), \ldots , \Gamma_{N+1} = 1_A(N+1))= p^k(1-p)^{N-k}. 
\end{equation}
This is the new basic belief assignment after conditioning on $H_s$ and can be called the {\em posterior} in this case.

We can now use this result to compute belief in certain events, but before we do that, we would like to discuss the chosen modeling. We defined the basic belief assignment $m$ without using that we know that someone in the population has $\Gamma$, namely the criminal $C$; this information is only added once we condition on $H_s$. Perhaps some readers might find this somewhat counterintuitive. Why not use the information that there is at least one individual in $X$ which has $\Gamma$ in the definition of the basic belief assignment? We now explain how this can be done, and show that it leads to the same belief assignment after proper conditioning.

Let $y_1, \ldots, y_{N+1}$ be such that $\sum_i y_i =k+1$, $k=0, \ldots, N$. We define the basic belief assignment 
$m'$ as follows:
\begin{eqnarray}
\label{basis3}
m'(C \in \{ i \;:\; y_i=1 \},  S=x &,& \Gamma_1 = y_1 , \ldots , \Gamma_{N+1}= y_{N+1}) =\\
& =&  \frac{1}{N+1} \frac{p^{k+1} (1-p)^{N-k}}{1-(1-p)^{N+1}}. \nonumber
\end{eqnarray}
This belief assignment gives us our belief in the characteristics $\Gamma_i$ conditioned on the event that at least one $y_i$ is equal to 1. It expresses no knowledge about the identity of $C$ other than that we know that $C$ has $\Gamma$. We denote the event in (\ref{basis3}) by $A'(x, y_1, \ldots, y_{N+1})$.

Next we condition on the event $H'_s=\{S=s, \Gamma_S=1\}$. As before, $H'_s \cap A'(x,y_1, \ldots, y_{N+1})$ can only be nonempty if $x=s$ and $y_s=1$. In that case $H'_s \cap A'(x,y_1, \ldots, y_{N+1})$ is exactly the event in 
(\ref{basis2}) and the correct normalization is 
\begin{equation}
\frac{1}{N+1}\sum_{j=0}^N {N \choose j } \frac{p^{j+1}(1-p)^{N-j}}{1-(1-p)^{N+1}} = \frac{p}{N+1}\cdot \frac{1}{1-(1-p)^{N+1}}.
\end{equation}
It now follows that after conditioning we obtain the same belief assignment as before, since the extra term $1-(1-p)^{N+1}$ appears in both the numerator and the denominator. We conclude that the two approaches lead to the same result, as they should.

Now that we have computed the new basic belief assignment, we can compute our belief in certain events, most notably our belief in the event that $C \in B$ for $B \subseteq X$ with $s \in B$. For this event we have:\footnote{We write $|B|$ for the number of elements in the set $B$.}
\begin{eqnarray}
\label{geleerd}
\B_{H_s}(C \in B) &=& \sum_{E \subseteq \{C \in B\}} m_{H_s}(E)  \nonumber \\
 & = & \sum_{A \subseteq B | s \in A}  m_{H_s}( C \in A ; S=s; \Gamma_1 = 1_A(1) \ldots ;\nonumber  \\
& & \ldots ; \Gamma_{N+1} = 1_A(N+1))\nonumber \\
& = &  \sum_{k=0}^{|B|-1} {|B|-1 \choose k} p^k (1-p)^{N-k} \nonumber \\
& = &  (1-p)^{N+1-|B|}.
\end{eqnarray}

An interesting special case occurs when $B=\{s\}$. The conditional belief that $C=s$ is apparently given by 
\begin{equation}
\label{ip1}
\B_{H_s}(C=s)=(1-p)^N,
\end{equation} 
simply take $|B|=1$. This formula has a simple interpretation: the belief that $s$ is the criminal is just the probability that all other members of the population are excluded since they have the wrong profile.

It is interesting to compare this answer to the classical one, in which a uniform prior is taken. In the classical case, the posterior probability that $C=s$ is equal to 
\begin{equation}
\label{ip2}
\frac{1}{1+Np},
\end{equation}
see e.g.\ \cite{sloomee11} or \cite{baldon95}. We observe that
\begin{equation}
\frac{1}{1+Np} = \prod_{k=0}^{N-1} \frac{1+kp}{1+(k+1)p} > \left(\frac{1}{1+p}\right)^N > \left(\frac{1-p^2}{1+p}\right)^N = (1-p)^N.
\end{equation}
Hence in our setting, the belief that $C=s$ is always smaller than in the classical case, something we can intuitively understand by recalling that we assign prior belief zero to this event. To give some indication of the difference between the two answers, if $p \sim N^{-1}$ (for $N \rightarrow \infty)$, then (\ref{ip2}) $\sim \frac{1}{2}$, while (\ref{ip1}) $\sim e^{-1}$.

Since belief functions generalize probability distributions, we should be able to re-derive the classical result (\ref{ip2}) using our approach, and we now show that this is indeed the case. If we want to take a uniform prior for the criminal, then the basic belief assignment, denoted by $m^c$, is as follows. Let $y_1, \ldots, y_{N+1}$ be such that $\sum_{i=1}^{N+1} y_i=k$. Then our prior is given by
\begin{equation}
\label{basis5}
m^c(C=t, S=x, \Gamma_1 = y_1 , \ldots , \Gamma_{N+1}= y_{N+1}) = \frac{1}{(N+1)^2} p^{k} (1-p)^{N+1-k},
\end{equation}
for all $t, x \in X$, and $m^c(A)=0$ for all other subsets of $\Omega$. Note that the corresponding belief function is a probability distribution, since only singletons have positive basic belief. Next we condition on the same event $H_s=\{S=s , \Gamma_C = \Gamma_S = 1\}$ as before. The intersection of the set in (\ref{basis4}) with $H_s$ is only non-empty if $s=x$ and $
y_s=y_t=1$. The probability that $s=x$ and $y_s=1$ is $p/(N+1)$ as before. Given this, the probability that also $y_t=1$ is 1 if $t=s$ and $p$ if $t\neq s$. Hence the intersection is non-empty with probability
$$
\frac{p}{N+1}\left( \frac{1}{N+1} + \frac{N}{N+1} p\right)= \frac{p(1+Np)}{(N+1)^2}.
$$
We can now compute the conditional belief assignment $m^c_{H_s}$ but we note that we only need (with $y_s=1$ and $\sum_i y_i=k+1$)
\begin{eqnarray*}
m^c_{H_s}(s, s, y_1, \ldots, y_{N+1}) &=& \frac{p^{k+1}(1-p)^{N-k}}{(N+1)^2} \frac{(N+1)^2}{p(1 + Np}\\
&=& \frac{p^{k}(1-p)^{N-k}}{1+Np}.
\end{eqnarray*}
Summing over all $k=0, \ldots, N$ and using Newton's binomium, we see that the conditional belief that $S=C=s$ is given by (\ref{ip2}), as required.

This example illustrates that we lose nothing by working with belief functions, and that belief functions only add flexibility. If a certain classical prior is reasonable then we can take that prior and work with it. If there are reasons to have a non-classical prior, for instance complete ignorance within a given population, then belief functions are flexible enough to deal with this.

\subsection{The cold case generalized}
\label{afh}
So far we have assumed that we can describe the realization of the characteristic with a classical probability distribution with independence between different individuals. The flexibility that is possible with belief functions was only used in order to model complete ignorance about the culprit. 

But it is possible to be even more flexible with our knowledge of the characteristics, and use belief functions also for them. This entails a certain uncertainty about the occurrence of the characteristic $\Gamma$. We write $p>0$ for the probability that we determine that an individual has the property. We write $q \geq 0$ for the probability that we determine that an individual does not have the property and write $r \geq 0$ for the probability that we can not determine if an individual has the property.  This then leads to the following prior basic belief assignment $m$ on $\Omega = X \times X \times \{0,1\}^{N+1}$ as the analogue of (\ref{basis}):
\begin{equation}
\label{basisnieuw}
m(S=x, \Gamma_1 \in Y_1 , \ldots , \Gamma_{N+1} \in Y_{N+1}) = \frac{1}{N+1} p^{k_1} q^{k_0} r^{N+1-k_0-k_1},
\end{equation}
where $Y_i \in \{ \{0,\}, \{1\}, \{0,1\}\}$ and $k_1 = \sum_i 1(Y_i = \{1\})$ and $k_0 = \sum_i 1(Y_i = \{0\})$.

From here on, the analysis is more or less as before. We denote the set in (\ref{basisnieuw}) by $A(x, Y_1, \ldots, Y_{N+1})$. 
Again we want to condition on $H_s = \{S =s \;:\; \Gamma_C = \Gamma_S =1 \}$. The intersection $H_s \cap A(x,Y_1,\ldots,Y_{N+1})$ is nonempty if and only if $x=s$ and $Y_s=\{1\}$ and can in that case be written as
\begin{equation}
\{ C \in \{i \;:\; Y_i \not= \{0\} \}, S=s, \Gamma_1 \in Y_1, \ldots, \Gamma_n \in Y_{N+1}\}.
\end{equation}
The total mass of these sets is
\begin{equation}
\frac{p}{N+1},
\end{equation}
for the same reason as before. Hence we have the posterior
\begin{eqnarray*}
m_{H_s}(A(s, Y_1, \ldots, Y_{N+1})\cap H_s) &=& \frac{p^{k} q^l r^{N+1-k-l}}{N+1} \cdot \frac{N+1}{p}\\
&=& p^{k-1} q^l r^{N+1-k-l}.
\end{eqnarray*}
For $s \in B$ we find
\begin{equation}
\B_{H_s}(C \in B) = q^{N+1-|B|}.
\end{equation}
When $r=0$ this is just (\ref{ip1}) as it should. It is noteworthy that the belief in $C \in B$ only depends on $|B|$ and $q$, and not on $p$ and $r$. This is to be expected: the belief that $s$ is the criminal is the probability that we know for {\em sure} that all other members of the population do not have the correct profile.

\subsection{The search case}

In the search variant, we do not choose a random individual $S$ but we check the inhabitants one by one in a random order, until an individual with the relevant characteristic is found. However, we only take into account the result of the search and not any information about the search itself. As a consequence, the search case boils down to picking a random individual from the subset of the population that has the characteristic. As in the cold case, there are (at least) two ways  to approach the situation. In the first approach, we do not immediately link the characteristic to the crime, and then it is a priori not certain that an individual with the characteristic found at the scene of the crime exists in the population. In the second approach, we condition the distribution of the characteristics on the fact that we know at least one individual has it. 

For the first approach, the space $\Omega$ must be rich enough to accommodate the possibility that no one has the characteristic. To this end, we set $\Omega = X \times (X \cup \{*\}) \times \{0,1\}^{N+1}$ and let $C: \Omega \rightarrow X,S: \Omega \rightarrow X \cup \{*\}$ and $\Gamma_i: \Omega \rightarrow \{0,1\}$ be projections on respectively the first, second and $i+2$-th coordinate. If no one has the characteristic we set $S=*$ to encode that we could not find a suspect (which is only the case if $\Gamma_1 = \Gamma_2 = \cdots = \Gamma_{N+1}= 0$). If at least one individual has the characteristic, it is clear that $S$ is uniformly selected from the subset of individuals that have the trait. Let $x \in X$. For $y_1,\ldots,y_{N+1}$ such that $\sum_i y_i=k$, $k=1, \ldots, N+1$, and $y_x=1$ we set our prior basic belief assignment as follows:
\begin{equation}
m(S=x, \Gamma_1 = y_1 , \ldots , \Gamma_{N+1}= y_{N+1}) = \frac{1}{k} p^{k} (1-p)^{N+1-k}
\end{equation}
and
$$
m(S=*, \Gamma_1 = \Gamma_2 = \ldots = \Gamma_{N+1} = 0)=(1-p)^{N+1}.
$$
Note that there are two differences compared to the cold case: we have to assume that $y_x=1$, and we have to divide by $k$ rather than by $N+1$.

Next we link the characteristic to the crime and we condition on $H_s=\{S=s , \Gamma_C = \Gamma_S = 1\}$. Note that the conditioning does not contain information about the length of the search or the identity of searched individuals. We only know that $s$ was the first one to be found with the characteristic. We have, for $y_1,\ldots,y_{N+1}$ such that $y_s=1$,
\begin{eqnarray}
\label{opnieuw}
& & \{ S=s; \Gamma_1 = y_1 ; \ldots ; \Gamma_{N+1}= y_{N+1} \} \cap H_s \\
& = & \{ C \in \{ i \;:\; y_i=1 \} , S=s; \Gamma_1 = y_1 ,\ldots , \Gamma_{N+1}= y_{N+1} \}. \nonumber
\end{eqnarray}
Note that the sets in (\ref{opnieuw}) are the only subsets of $\Omega$ with positive mass that have a nonempty intersection with $H_s$. Hence the normalization follows from the total mass of such sets, which is equal to
\begin{equation}
\label{berekening}
\begin{aligned}
M_s & := \sum_{j=0}^N {N \choose j }\frac{1}{j+1} p^{j+1}(1-p)^{N-j} \\
& = \sum_{j=0}^N \frac{1}{N+1} {N+1 \choose j+1 } p^{j+1}(1-p)^{N-j} \\
& = \frac{1- (1-p)^{N+1}}{N+1}.
\end{aligned}
\end{equation}
Notice that we can also derive (\ref{berekening}) by observing that $M_s$ does not depend on $s$ and thus $(N+1)M_s$  plus $\B(S=*)=(1-p)^{N+1}$ adds up to $1$.

Hence for $A \subseteq X$ such that $|A|=k+1$ and $s \in A$ we have the posterior
\begin{eqnarray}
\label{weert}
& & m_{H_s}( C \in A , S=s, \Gamma_1 = 1_A(1) ; \ldots , \Gamma_{N+1} = 1_A(N+1)) \nonumber \\
&= &\frac{ p^{k+1} (1-p)^{N-k}}{k+1} \frac{N+1}{ 1- (1-p)^{N+1}}.
\end{eqnarray}
Now we can compute our belief in certain events. For instance, our belief in $\{C=s\}$ is given by
\begin{eqnarray}
\label{antwoordd}
\B_{H_s}(C=s) &=& m_{H_s}( C = s, S=s, \Gamma_s=1,  \cap_{i \not=s} \{ \Gamma_i = 0 \})\nonumber \\
&=&  \frac{p (1-p)^{N}(N+1)}{ 1- (1-p)^{N+1}},
\end{eqnarray}
which is different compared to the cold case. There is a very natural interpretation of the expression in (\ref{antwoordd}). In the numerator, we have the probability that a binomially distributed random variable with parameters $N+1$ and $p$ is equal to 1. The denominator is the probability that this random variable is positive, so (\ref{antwoordd}) is the conditional probability for such a random variable to be 1 given it is positive. This makes sense, since we can only know for sure that $C=s$ when $s$ is the only one with the characteristic. Notice that (\ref{antwoordd}) equals zero if $p=1$ and if $p<1$ we can rewrite (\ref{antwoordd}) as
\begin{equation}
\label{antwoorddharmonic}
\frac{N+1}{\sum_{k=0}^N (1-p)^{-k}}.
\end{equation}

As already mentioned above, there is the alternative approach in which we deduce from the characteristic found at the crime scene that at least one individual has the characteristic. In this case we do not need the extended $\Omega$ since $S=*$ has probability zero. For $y_1, \ldots, y_{N+1}$ such that $y_x=1$ and $\sum_i y_i= k+1$, $k=0, \ldots, N$, the basic belief assignment is now given by
\begin{eqnarray*}
\label{basis4}
m'(C \in \{ i \;:\; y_i=1 \},  S=x &,& \Gamma_1 = y_1 , \ldots , \Gamma_{N+1}= y_{N+1}) =\\
& =&  \frac{1}{k+1} \frac{p^{k+1} (1-p)^{N-k}}{1-(1-p)^{N+1}}.
\end{eqnarray*}
After this, we condition on $H'_s=\{S=s\}$.  It is easily seen that the normalizing factor is just $1/(N+1)$, and this immediately leads to the same formula as in (\ref{antwoordd}).

In the classical case, starting with a uniform probability distribution, the posterior probability that $C=s$ is equal to
\begin{equation}
\label{klassiek2}
\frac{1-(1-p)^{N+1}}{(N+1)p} = \frac{1}{N+1} \sum_{j=0}^N (1-p)^j,
\end{equation}
see e.g.\ \cite{baldon95}. Note that (\ref{klassiek2}) is the arithmetic mean of $1,1-p,(1-p)^2, \ldots , (1-p)^N$, while (\ref{antwoorddharmonic}) is the harmonic mean of the same sequence. Hence the answer using our approach is - as it was in the cold case - smaller than the classical answer.

We briefly demonstrate that we can also derive this classical result with our technology. In the classical case, the basic belief assignment is given by
\begin{equation}
\label{basis7}
m^c(C=c, S=x, \Gamma_1=y_1, \ldots, \Gamma_{N+1}=y_{N=1})=\frac{1}{k} \frac{1}{N+1} p^k (1-p)^{N+1-k},
\end{equation}
whenever $\sum_i y_i =k$, and 
\begin{equation}
m^c(S=*, \Gamma_1 = \Gamma_2 = \cdots = \Gamma_{N+1} = 0)=(1-p)^{N+1}.
\end{equation}
Note that the corresponding belief function is a probability measure, since only singletons have positive basic belief assignments. We condition on $H_s=\{S=s, \Gamma_C =\Gamma_s=1\}$ as usual. To compute the (classical) belief in $C=s$ conditioned on $H_s$ we need to compute the conditional basic belief assignment $m^c_{H_s}$, for which we need the correct normalizing constant and the total mass of sets in (\ref{basis7}) whose intersection with $H_s$ is not empty. Elementary combinatorics gives that the latter is equal to
\begin{equation}
\label{uitkomst}
\sum_{k=0}^N {N \choose k} \frac{1}{k+1} \frac{1}{N+1} p^{k+1} (1-p)^{N-k}= \frac{1-(1-p)^{N+1}}{(N+1)^2}. 
\end{equation}
The normalizing constant follows from the total mass assigned to sets whose intersection with $H_s$ is non-empty. For this intersection to be non-empty, we need that $x=s$, $y_s=1$ and $y_C=1$. In the classical probabilistic experiment described by $m^c$, we simply need to compute the probability that $H_s$ occurs. We need subsequently (1) $y_C=1$, and (2) $s$ is chosen (implying that $y_s=1$). The first step occurs with probability $p$. Given this, we now know that not all labels are 0, and every individual has the same probability to be chosen. Hence, the conditional probability of Step (2) given Step (1) is simply $1/(N+1)$. It follows that the correct normalizing constant is $p/(N+1)$. Combining this with (\ref{uitkomst}) yields (\ref{klassiek2}).

\section{Parental identification}
\label{parental}

Our next example concerns the situation in which we have a known mother and a known child, but we do not know who the father is. We assume that there is a set $X=\{1,\ldots,N+1\}$ of potential fathers. We would like to make belief statements about the possible father-ship of someone chosen from the population, based on the DNA profile of this chosen person. In order to keep things as simple as possible, we assume that we only consider one specific locus of the DNA. Furthermore, we assume that the alleles of mother and child at that locus are such that we know what the paternal allele must be. Every potential father $i$ in $X$ has two alleles at the locus, and we denote by $\Gamma_i$ the number of `correct' alleles matching the paternal allele, hence $\Gamma_i \in \{0,1,2\}.$   

In order to set up our prior belief function, we set $\Omega = X \times X \times \{0,1\} \times \{0,1,2\}^{N+1}$ and let $F: \Omega \rightarrow X,S: \Omega \rightarrow X$, $A: \Omega \rightarrow \{0,1\}$ and $\Gamma_i: \Omega \rightarrow \{0,1,2\}$ be projections on respectively the first, second, third and $i+3$-th coordinate. $F$ represents the father, and $S$ the selected individual which is the putative father. The indicator $A$ is associated with the putative father, and is only relevant if the putative father has exactly one `correct' allele; in that case $A=1$ means that he passes on this correct allele to his child, while $A=0$ indicates that he does not do that. Let $p_0,p_1$ and $p_2$ be the probabilities that an individual has respectively $0,1,2$ alleles of the right type. We assume that we know these probabilities from population surveys. Set $k_j := \sum_{i=1}^{N+1} 1(y_i=j)$ 
for $j=0,1,2$. Then the prior basic belief assignment is given by
\begin{equation}
\label{basic8}
m(S=x,  A=a, \Gamma_1 = y_1 , \ldots , \Gamma_{N+1}= y_{N+1}) = \frac{1}{2}\frac{1}{N+1} p_0^{k_0} p_1^{k_1} p_2^{k_2}.
\end{equation}
As before, this basic belief assignment is a summary of what we know, and is formulated in terms of items that can be well described by classical probabilities. The factor $\frac12$ comes from the fact that a father passes a randomly chosen allele to his child. Note that the basic belief assignment in (\ref{basic8}) only assigns positive basic belief to sets that contain no information whatsoever about the father, other than the fact that he belongs to the given population. We denote the event in (\ref{basic8}) by $E(x,a, y_1, \ldots, y_{N+1})$.

Now we have to distinguish between two scenarios: the suspect in question has one or two alleles of the right type. We look at the case he has two such alleles first. This means we want to condition on
\begin{equation}
H_{s,2}=\{S=s , \Gamma_S = 2, \{\Gamma_F=2\} \cup \; \{ \Gamma_F=1, A=1 \}\}.
\end{equation}
We have to assign the mass originally assigned to $E(x,a, y_1, \ldots, y_{N+1})$ to the intersection of this set with $H_{s,2}$, and normalize suitably. This intersection is non-empty precisely when $x=s$ and $y_s=2$, and as before, the total mass of this can be calculated by viewing the basic belief assignment in (\ref{basic8}) as a classical probabilistic experiment, leading to a normalizing constant of
\begin{equation}
\frac{p_2}{N+1}.
\end{equation}
It follows that if $y_s=2$ we have the posterior
\begin{eqnarray*}
m_{H_{s,2}}(E(s,a, y_1, \ldots, y_{N+1})\cap H_{s,2}) &=& \frac{1}{2}\frac{1}{N+1} p_0^{k_0} p_1^{k_1} p_2^{k_2} \cdot \frac{N+1}{p_2}\\
&=& \frac{1}{2}  p_0^{k_0} p_1^{k_1} p_2^{k_2-1}.
\end{eqnarray*}
This basic belief assignment leads to the following posterior belief function. Let $B \subseteq X$ be such that $s \in B$. If $a=1$, then we are back in the classical cold case of the island problem, with $p$ replaced by $1-p_0$, since one correct allele will be enough now:
\begin{equation}
\begin{aligned}
& \B_{H_{s,2}}(F \in B, A=1) \\
& = \frac{1}{2} \sum_{k_0=0}^{|B|-1} {|B| - 1 \choose k_0} p_0^{k_0} \sum_{k_1=0}^{|B|-1-k_0} {|B| - 1 - k_0 \choose k_1} p_1^{k_1} p_2^{N-k_1-k_0} \\
& = \frac{1}{2} \sum_{k_0=0}^{|B|-1} {|B| - 1 \choose k_0} p_0^{k_0} (p_1+p_2)^{N-k_0} \\
& = \frac{1}{2} p_0^{N+1-|B|}.
\end{aligned}
\end{equation}
If $a=0$ then the situation reduces to the cold case of the island problem with $p=p_2$, since a potential father now needs two alleles. By an analogous computation we obtain:
\begin{equation}
\B_{H_{s,2}}(F \in B, A=0) =  \frac{1}{2} (p_0+p_1)^{N+1-|B|}.
\end{equation}
Hence
\begin{equation}
\label{belfun1}
\begin{aligned}
\B_{H_{s,2}}(F \in B) & = \B_{H_{s,2}}(F \in B, A=1) + \B_{H_{s,2}}(F \in B, A=0) \\
& = \frac{1}{2} \left( p_0^{N+1-|B|} + (p_0+p_1)^{N+1-|B|} \right).
\end{aligned}
\end{equation}
In the special case that $B=\{s\}$, we find
\begin{equation}
\label{antwoord2}
\B_{H_{s,2}}(F=s) = \frac{1}{2} p_0^{N} + \frac{1}{2}(p_0+p_1)^{N},
\end{equation}

Next we consider the case in which we condition on
\begin{equation}
H_{s,1}=\{S=s , \Gamma_S = 1, \{\Gamma_F=2\} \cup \; \{ \Gamma_F=1; A=1 \}\}.
\end{equation}
Now $E(x,a, y_1, \ldots, y_{N+1}) \cap H_{s,1}$ can only be non-empty if $x=s$ and $y_s=1$, which happens with probability $\frac{p_1}{N+1}$. Note however that, among the intersections, if $a=0$ and $y_i \not=2$ for all $i$ then the intersection is in fact empty. The latter occurs with probability $\frac12(p_0+ p_1)^N$, and hence the correct normalizing constant is
\begin{equation}
\frac{p_1}{N+1} - \frac{1}{2} \frac{p_1}{N+1} (p_0+p_1)^{N} = \frac{p_1}{N+1} \left( 1 - \frac{1}{2}(p_0+p_1)^{N} \right).
\end{equation}
It follows (except when $a=0$ and $y_i \neq 2$ for all $i$) that we have the posterior
\begin{equation}
m_{H_{s,1}}(E(s,a, y_1, \ldots, y_{N+1})\cap H_{s,1}) 
= \frac{1}{2} \frac{ p_0^{k_0} p_1^{k_1-1} p_2^{k_2}}{ 1 - \frac{1}{2}(p_0+p_1)^{N} }.
\end{equation}
We can now compute the corresponding posterior belief function using a similar argument as the one which led to 
(\ref{belfun1}). When $a=1$ we have the same situation as before. When $a=0$, we only need to rule out the case in which all labels $y_i$ are 0 or 1. This means that compared to the combinatorics of the previous case, we have to subtract $(p_0+p_1)^N$ in this case. For $s \in B$ we now get
\begin{equation}
\B_{H_{s,1}}(F \in B) = \frac{1}{2} \frac{ p_0^{N+1-|B|} + (p_0+p_1)^{N+1-|B|} - (p_0+p_1)^N}{1- \frac{1}{2}(p_0+p_1)^N}.
\end{equation}
In the special case in which $B=\{s\}$ we find
\begin{equation}
\label{antwoord1}
\B_{H_{s,1}}(F=s) =  \frac{ \frac{1}{2} p_0^{N} }{1-\frac{1}{2}(p_0+p_1)^N}
\end{equation}
and we note that (\ref{antwoord1}) is smaller than (\ref{antwoord2}), which is what common sense requires.

Finally, we compare our answers with the answers we get when we start with a uniform prior instead of a uninformative prior. Since we could not find our specific situation in the literature, we perform the necessary calculations. To that end, we set the basic belief assignment
\begin{equation}
\label{basic8klas}
\begin{aligned}
& m^c(S=x, F=y,  A=a, \Gamma_1 = y_1 , \ldots , \Gamma_{N+1}= y_{N+1}) \\
& = \frac{1}{2}\frac{1}{(N+1)^2} p_0^{k_0} p_1^{k_1} p_2^{k_2},
\end{aligned}
\end{equation}
which only assigns positive support to singletons and thus corresponds to a probability distribution. By computing
\begin{equation}
\B^c(\{F=s\} \cap H_{s,1}) = \frac{1}{(N+1)^2}\frac{1}{2}p_1
\end{equation}
and
\begin{equation}
\begin{aligned}
\B^c(H_{s,1}) & = \B^c(H_{s,1} \cap \{F=s\}) + \B^c(H_{s,1} \cap \{F \not= s\}) \\
& =  \frac{1}{(N+1)^2}p_1 \frac{1}{2} +  \frac{1}{(N+1)^2} p_1  N\left(\frac{1}{2}p_1+p_2 \right),
\end{aligned}
\end{equation}
we find that
\begin{equation}
\begin{aligned}
\B^c_{H_{s,1}}(F=s) & = \frac{\B^c(\{F=s\} \cap H_{s,1})}{\B^c(H_{s,1})} \\
& = \frac{1}{1+ N(p_1 + 2p_2)} \\
& \geq \B_{H_{s,1}}(F=s).
\end{aligned}
\end{equation}
Similarly, we get
\begin{equation}
\begin{aligned}
\B^c_{H_{s,2}}(F=s) & = \frac{\B^c(\{F=s\} \cap H_{s,2})}{\B^c(H_{s,2})} \\
& = \frac{1}{1+ N(\frac{1}{2}p_1 + p_2)} \\
& \geq \B_{H_{s,2}}(F=s).
\end{aligned}
\end{equation}

\section{An analogue of Bayes' rule}
\label{rule}
In the context of belief functions, we do not have a useful explicit closed expression as the analogue of Bayes' rule as formulated  in (\ref{bayesintro}). However, we do have a similar process of updating prior belief statements when we are given additional knowledge. In this section we explain how the analogue of working with Bayes' rule works in our context, and how we can use it in computations with belief functions. We think this is best explained by analyzing a number of instructive example, and this is what we turn to now.

We recall the odds form of Bayes' rule:
\begin{equation}
\label{bayesjawel}
\frac{P(G|E)}{P(\bar{G}|E)}= \frac{P(E|G)}{P(E|\bar{G})} \cdot \frac{P(G)}{P(\bar{G})},
\end{equation}
where $\bar{G}$ denotes the negation or complement of $G$. (The rule can also be used when $\bar{G}$ is any other set.)  Whatever the interpretation of the classical probability measure $P$, subjective, frequentistic or otherwise, $P(G|E)$ represents the posterior probability of guilt conditioned on the available evidence, while $P(G)$ denotes the prior probability of guilt, before taking the evidence $E$ into account. 

One can, however, also look at the posterior odds from the perspective of the joint distribution of $G$ and $E$. Once we have the full joint distribution, the posterior odds (and any other quantity relating $E$ and $G$, for that matter) can be computed. Equation (\ref{bayesjawel}) in fact tells us that in the classical setup, the complete joint distribution of $G$ and $E$ can only be determined as a joint effort of both expert (likelihood ratio) and legal representative (prior). They both contribute to the joint distribution, albeit in different ways, and as such both contribute to the posterior odds as well. This means that to develop the analogue of Bayes' rule in our context of belief functions, we need to explain how the joint belief of $G$ and $E$ can be constructed, under a given prior belief on $G$ (which may be complete ignorance). 

Both $G$ and $E$ can, in our examples, only take values 0 and 1, where a 1 represents guilt and the existence of evidence, respectively. So we work on $\Omega= \{0,1\} \times \{0,1\}$ where, for instance, the state $(1,0)$ corresponds to the situation that suspect is guilty, but there is no evidence against him or her. We start with a classical situation which we first compute with (\ref{bayesjawel}) and after that a second time from the perspective of belief functions. After that we present two examples which cannot be done in a classical way. 

\begin{example} ({\em A classical situation})
\label{exclassic}
Suppose we would like to impose a prior probability of $0.1$ that our suspect is guilty, and a prior probability $0.9$ that he is not. We assume further that if the suspect is guilty, evidence is found with probability $0.8$, and if he is not guilty, evidence is found with probability $0.05$. This leads to prior odds of $1/9$ and a likelihood ratio of $0.8/0.05=16$. Hence the posterior odds are 16/9, and this corresponds to a posterior probability of guilt of $16/25=0.64$. 

Now let us recompute this example in the context of belief functions. In this classical case this may appear as more work, but the point is of course that the forthcoming method can be extended to non-classical situations.

The assumptions imply that our basic belief assignment on $\{0,1\}^2$ should be such that when we project on the first coordinate, we must have total mass $0.1$ on $G=1$ and mass $0.9$ on $G=0$. In combination with the other assumptions about the probabilities to find evidence under various hypotheses, this leads to the following basic belief assignment on $\{0,1\}^2$:

\begin{equation}
\begin{aligned}
m \left( \; \begin{tabular}{ | l | c | c | c |}
    \hline
 & $G=0$ & $G=1$ \\ \hline
$E=0$ &  $*$   &     \\ \hline
$E=1$ &   &      \\ \hline
\end{tabular} \; \right) = 0.9 \times 0.95=0.855,
\end{aligned}
\end{equation}

\begin{equation}
\begin{aligned}
m \left( \; \begin{tabular}{ | l | c | c | c |}
    \hline
 & $G=0$ & $G=1$ \\ \hline
$E=0$ &     &     \\ \hline
$E=1$ & $*$  &       \\ \hline
\end{tabular} \; \right) = 0.9 \times 0.05=0.045,
\end{aligned}
\end{equation}

\begin{equation}
\begin{aligned}
m \left( \; \begin{tabular}{ | l | c | c | c |}
    \hline
 & $G=0$ & $G=1$ \\ \hline
$E=0$ &     &   $*$  \\ \hline
$E=1$ &   &      \\ \hline
\end{tabular} \; \right) = 0.1 \times 0.2 = 0.02,
\end{aligned}
\end{equation}

\begin{equation}
\begin{aligned}
m \left( \; \begin{tabular}{ | l | c | c | c |}
    \hline
 & $G=0$ & $G=1$ \\ \hline
$E=0$ &     &     \\ \hline
$E=1$ &   &   $*$    \\ \hline
\end{tabular} \; \right) = 0.1 \times 0.8 = 0.08.
\end{aligned}
\end{equation}
Conditioning on $E=1$ means we disregard the first and third option since these sets have an empty intersection with $E=1$. Only the last option gives mass to $G=1$ after conditioning, and therefore the conditional probability that $G=1$ given $E=1$ is given by $0.08/(0.08 + 0.045)= 0.64$, agreeing with the earlier computation.
\end{example}

Suppose now that we want to model complete ignorance about the guilt of our suspect. This boils down to the requirement that our basic belief assignment on $\Omega$ must be such that when we project any set with positive basic belief assignment on the first coordinate, we must get the full marginal $\{0,1\}$. There are nine subsets of $\Omega$ with this property, and it is very instructive to list them and to decide whether or not each of them should be assigned positive basic belief assignment. This is more work than assigning conditional probabilities for application of the classical Bayes' rule, but the advantage is that we can very precisely model the possibly subtle relation between $G$ and $E$. We list the possibilities as 2 by 2 arrays, hopefully increasing readability and clarity of the exposition. 

We start with the three `rectangular' sets, whose interpretation is straightforward:

\medskip
I:
\begin{tabular}{ | l | c | c |}
    \hline
 & $G=0$ & $G=1$ \\ \hline
$E=0$ & $*$ & $*$  \\ \hline
$E=1$ & $*$ & $*$  \\ \hline
\end{tabular}

\medskip
Interpretation: uninformative.

\medskip\noindent

\medskip\medskip
II:
\begin{tabular}{ | l | c | c |}
    \hline
 & $G=0$ & $G=1$ \\ \hline
$E=0$ & $*$ & $*$  \\ \hline
$E=1$ &  &   \\ \hline
\end{tabular}

\medskip
Interpretation: absence of evidence in all circumstances.

\medskip\noindent

\medskip\medskip
III:
\begin{tabular}{ | l | c | c |}
    \hline
 & $G=0$ & $G=1$ \\ \hline
$E=0$ &  &   \\ \hline
$E=1$ & $*$ & $*$  \\ \hline
\end{tabular}

\medskip
Interpretation: existence of evidence in all circumstances. 

\medskip\medskip
The next possibility corresponds to the situation that evidence exists precisely when suspect is guilty. In ideal cases one would like to assign a lot of belief to this set. 

\medskip
IV:
\begin{tabular}{ | l | c | c |}
    \hline
 & $G=0$ & $G=1$ \\ \hline
$E=0$ & $*$ &   \\ \hline
$E=1$ &  & $*$  \\ \hline
\end{tabular}

\medskip
Interpretation: incriminating evidence that exists precisely when suspect is guilty.

\medskip\medskip
There are two perturbations of Option IV, namely with the possibilities of a false positive respectively false negative added:

\medskip\medskip
V:
\begin{tabular}{ | l | c | c |}
    \hline
 & $G=0$ & $G=1$ \\ \hline
$E=0$ & $*$ &  \\ \hline
$E=1$ & $*$ & $*$  \\ \hline
\end{tabular}

\medskip
Interpretation: as in Option IV but with false positive for incriminating evidence possible. Alternatively, we can interpret this as ruling out a false negative for incriminating evidence.

\medskip\medskip
VI:
\begin{tabular}{ | l | c | c |}
    \hline
 & $G=0$ & $G=1$ \\ \hline
$E=0$ & $*$ & $*$  \\ \hline
$E=1$ & & $*$  \\ \hline
\end{tabular}

\medskip
Interpretation: as in Option IV but with false negative for incriminating evidence possible. Alternatively, we can interpret this as ruling out a false positive for incriminating evidence.

\medskip\medskip
The remaining options seem more suitable for exculpatory evidence:

\medskip\medskip
VII:
\begin{tabular}{ | l | c | c |}
    \hline
 & $G=0$ & $G=1$ \\ \hline
$E=0$ &  & $*$  \\ \hline
$E=1$ & $*$ &   \\ \hline
\end{tabular}

\medskip
Interpretation: perfect exculpatory evidence.

\medskip\medskip
The last two options are perturbations of Option VII.
 
\medskip\medskip
VIII:
\begin{tabular}{ | l | c | c |}
    \hline
 & $G=0$ & $G=1$ \\ \hline
$E=0$ &  & $*$  \\ \hline
$E=1$ & $*$ & $*$  \\ \hline
\end{tabular}

\medskip
Interpretation: as Option VII with false positive for exculpatory evidence possible. Alternatively, we can interpret this as ruling out a false negative for exculpatory evidence.

\medskip\medskip
IX:
\begin{tabular}{ | l | c | c |}
    \hline
 & $G=0$ & $G=1$ \\ \hline
$E=0$ & $*$ & $*$  \\ \hline
$E=1$ & $*$ &   \\ \hline
\end{tabular}

\medskip
Interpretation: as Option VII with false negative for exculpatory evidence possible.  Alternatively, we can interpret this as ruling out a false positive for exculpatory evidence.

\medskip\medskip
Obviously, the actual belief that one wants to assign to various possibilities depends on personal taste, background knowledge, and on details of the situation. This is not different from the situation when working in a classical setting. In the context of our belief functions, we need to decide about the assignment of prior belief to a relatively large number of sets. In our example above, we need to make decisions about 9 sets. To compare: in our example, a classical Bayesian analysis requires no more than 4 probabilities, namely the probabilities that $E=i$ if $G=j$, for $i,j \in \{0,1\}$. But our belief functions are much more flexible, and allow us to tailor the modeling very much to the actual situation. 

Once we have decided about our basic belief assignment, we can perform computations with the calculus that we have developed in Section \ref{twee}. 

\begin{example} ({\em The cold case revisited})
\label{exignorance}
Consider the cold case from Section \ref{cold} once again. Suppose $s$ is chosen, and $\{G=1\}$ corresponds to $s$ being guilty. The event $\{E=1\}$ corresponds to the event that both $s$ and $C$ have the characteristic. We can now write down our basic belief assignment as follows.

First of all, 

\begin{equation}
\begin{aligned}
m \left( \; \begin{tabular}{ | l | c | c | c |}
    \hline
 & $G=0$ & $G=1$ \\ \hline
$E=0$ &  $*$   &   $*$  \\ \hline
$E=1$ &   &       \\ \hline
\end{tabular} \; \right) = 1-p,
\end{aligned}
\end{equation}
since this is our belief in the set that there is no evidence for sure, which is precisely the case when $s$ does not have the characteristic. Then we have

\begin{equation}
\begin{aligned}
m \left( \; \begin{tabular}{ | l | c | c | c |}
    \hline
 & $G=0$ & $G=1$ \\ \hline
$E=0$ &     &     \\ \hline
$E=1$ &  $*$ &  $*$     \\ \hline
\end{tabular} \; \right) = p^{N+1},
\end{aligned}
\end{equation}
since we know that there will be incriminating evidence precisely when everyone has the characteristic. Next we set 

\begin{equation}
\begin{aligned}
m \left( \; \begin{tabular}{ | l | c | c | c |}
    \hline
 & $G=0$ & $G=1$ \\ \hline
$E=0$ &  $*$   &   $*$  \\ \hline
$E=1$ &   &   $*$    \\ \hline
\end{tabular} \; \right) = p(1-p)^N,
\end{aligned}
\end{equation}
because this option corresponds to ruling out a false positive, which occurs when $s$ is the only one with the characteristic. The remaining mass goes to

\begin{equation}
\begin{aligned}
m \left( \; \begin{tabular}{ | l | c | c | c |}
    \hline
 & $G=0$ & $G=1$ \\ \hline
$E=0$ &  $*$   &   $*$  \\ \hline
$E=1$ &  $*$ &   $*$    \\ \hline
\end{tabular} \; \right) = p-p^{N+1} -p(1-p)^N,
\end{aligned}
\end{equation}
which makes perfect since, since this is the probability that $s$ has the characteristic, and in addition there is at least one other person with the characteristic, and also at least one who does not have it. In this situation we have complete ignorance when it comes to $s$. (Note that for $N=1$, the mass of the last option is 0. This is understandable, since we need at least three individuals for this option to be realized.)

Next we condition on $E=1$. The first set has empty intersection with $\{E=1\}$, and following the rules of our calculus we compute that the conditional belief that $G=1$ is equal to $(1-p)^N$, in accordance with (\ref{ip1}).

\end{example}

What we have done here is analogous to the classical situation. We started out with a belief function on $\{0,1\}^2$ which represented complete ignorance about the question whether or not the suspect is guilty, and which assigned positive basic belief to various possibilities. This is the analogue of a prior probability in the classical situation, together with the value of various conditional probabilities. Next we have conditioned on the event 
that $E=1$ and we have calculated the ensuing conditional belief in $G=1$. This is the analogue of calculating the conditional probability that $G=1$ given $E=1$ using Bayes' rule. Although we do not have an explicit formula like Bayes' rule at our disposal, we can nevertheless carry out the necessary calculations, leading to a well-defined answer. 

We would like to emphasize once more that we lose nothing when we work with belief functions, compared to the classical probabilistic set-up. We can start out with a classical probability distribution as our prior, and mimic a classical situation. We illustrate this with an example.

The relevant sets in Example \ref{exignorance} were completely disjoint from the relevant sets in Example \ref{exclassic}, but it may very well happen that we want to model both partial prior ignorance together with some partial specific prior belief. To this end we give our third example.
    
\begin{example} ({\em A mixture of ignorance and specific belief})
Consider the same situation as in the previous examples, but this time there is a prior indication that suspect is not guilty, say by a witness statement. Suppose we want to assign a total belief of $0.7$ to this witness. Then the basic belief assigment $m$ should assign a total mass of $0.7$ to sets whose projection onto $G$ is $\{0\}$. We have no indication about guilt, so $m$ should assign the remaining mass of $0.3$ to sets whose projection onto $G$ is $\{0,1\}$.

We now need to model what the probability of evidence is when suspect is not guilty. Suppose we have some reason to believe in false positives, but also some reason not to have any information about $E$ at all when $G=0$. These considerations could for example lead to the following prior belief assignment: 

\begin{equation}
\begin{aligned}
m \left( \; \begin{tabular}{ | l | c | c | c |}
    \hline
 & $G=0$ & $G=1$ \\ \hline
$E=0$ &  $*$   &    \\ \hline
$E=1$ &   &       \\ \hline
\end{tabular} \; \right) = 0.5,
\end{aligned}
\end{equation}

\begin{equation}
\begin{aligned}
m \left( \; \begin{tabular}{ | l | c | c | c |}
    \hline
 & $G=0$ & $G=1$ \\ \hline
$E=0$ &    &     \\ \hline
$E=1$ &  $*$ &       \\ \hline
\end{tabular} \; \right) = 0.05,
\end{aligned}
\end{equation}

\begin{equation}
\begin{aligned}
m \left( \; \begin{tabular}{ | l | c | c | c |}
    \hline
 & $G=0$ & $G=1$ \\ \hline
$E=0$ &  $*$   &     \\ \hline
$E=1$ &  $*$ &       \\ \hline
\end{tabular} \; \right) = 0.15,
\end{aligned}
\end{equation}

\begin{equation}
\begin{aligned}
m \left( \; \begin{tabular}{ | l | c | c | c |}
    \hline
 & $G=0$ & $G=1$ \\ \hline
$E=0$ &  $*$   &     \\ \hline
$E=1$ &   &   $*$    \\ \hline
\end{tabular} \; \right) = 0.2,
\end{aligned}
\end{equation}

\begin{equation}
\begin{aligned}
m \left( \; \begin{tabular}{ | l | c | c | c |}
    \hline
 & $G=0$ & $G=1$ \\ \hline
$E=0$ &  $*$   &     \\ \hline
$E=1$ &  $*$ &  $*$     \\ \hline
\end{tabular} \; \right) = 0.05,
\end{aligned}
\end{equation}

\begin{equation}
\begin{aligned}
m \left( \; \begin{tabular}{ | l | c | c | c |}
    \hline
 & $G=0$ & $G=1$ \\ \hline
$E=0$ &  $*$   &   $*$  \\ \hline
$E=1$ &   &   $*$    \\ \hline
\end{tabular} \; \right) = 0.05.
\end{aligned}
\end{equation}

When we now condition on $E=1$, the mass 0.5 of the first option is lost, so we have to renormalize by multiplying with 2. Hence we find that 
$$
m_{\{E=1\}}(G=0)=2(0.05 + 0.15) = 0.4,
$$
$$
m_{\{E=1\}}(G=1)=2(0.2 + 0.05) = 0.5,
$$
and
$$
m_{\{E=1\}}(G\in \{0,1\})=2 \cdot 0.05 = 0.1.
$$
The original belief of 0.7 not being guilty by in the witness statement is reduced to 0.4 after the evidence, but note that uncertainty persists. 
\end{example}    

Working with belief functions rather than probability distributions may seem rather complicated. Due to the nature of the theory in 
which we assign probabilities to subsets rather than to singletons, there are 
many more objects which can receive positive mass in our set-up compared to 
the classical one. However, in many cases modeling assumptions can be made 
which can, and should for that matter, drastically reduce the number of sets 
with positive basic belief assignment. The  computations themselves can of course be programmed.

Having explained how the analogue of applying Bayes' rule works in our context, a next step would be to investigate the analogue of a Bayesian network. Since we know how to mimic the working of Bayes' rule, we can now also create the analogue of a Bayesian network. This however, needs more space than we have in the current paper, and we carry out the construction of such a Bayesian networks in our context in a forthcoming paper.  

\section{Discussion and legal practice}
\label{practice}

How does the theory of belief functions as developed in this article relate to the legal practice when expert witnesses are called to testify? It is to this question that we now turn. 

We start by describing the current practice when classical probability theory is used. In this current practice, the legal representative (judge or jury member for instance) and expert witness play a very different role, and their contributions are well-separated. This is theoretically backed up by Bayes' rule. 
To be precise, writing $G$ for the event, say, that a certain suspect is the donor of a DNA profile found at the scene of the crime), and $E$ for the available evidence, the so called odds form of Bayes' rule states that
\begin{equation}
\label{bayes}
\frac{P(G|E)}{P(\bar{G}|E)}= \frac{P(E|G)}{P(E|\bar{G})} \cdot \frac{P(G)}{P(\bar{G})},
\end{equation}
where $\bar{G}$ denotes the negation or complement of $G$. Whatever the interpretation of the classical probability measure $P$, subjective, frequentistic or otherwise, $P(G|E)$ represents the posterior probability of $G$ conditioned on the available evidence, while $P(G)$ denotes the prior probability of $G$, before taking the evidence $E$ into account. 

Formula (\ref{bayes}) describes how the prior odds $P(G)/P(\bar{G})$ are transformed into the posterior odds
$P(G|E)/P(\bar{G}|E)$, by multiplying the prior odds with the so-called likelihood ratio $P(E|G)/P(E|\bar{G})$. 
This likelihood ratio is in the provenance of the forensic expert. This expert can, in certain circumstances at least, compute, estimate or assign a probability to the evidence under various hypothesis, and thereby compute the likelihood ratio. The expert does not express how likely these hypotheses themselves are, that is, the expert says nothing about the prior probabilities $P(G)$ and $P(\bar{G})$. Hence, the expert does not make a statement about the posterior probabilities $P(G|E)$ and $P(\bar{G}|E)$. Indeed, for the latter probabilities, one in addition needs the priors and these are not in the provenance of the expert. 

One can also look at the posterior odds from the perspective of the joint distribution of $G$ and $E$. Once we have the full joint distribution, the posterior odds (and any other quantity relating $E$ and $G$, for that matter) can be computed. But it is only through the joint effort of expert and legal representative that the joint distribution can be computed. In other words, the complete joint distribution of $G$ and $E$ can only be determined as a joint effort of both expert and legal representative. They both contribute to the joint distribution, albeit in different ways, and as such both contribute to the posterior odds as well.

Obviously, the legal practice is more complicated than this simple procedure suggests. Although Bayes' rule is mathematically not difficult, confusion between likelihood ratio and posterior odds arises easily, and there are many examples in which such confusion has had serious consequences see e.g.\ \cite{dawid} and the references therein. But also, as already articulated in the introduction, there are many situations possible where a prior in the context of classical probability theory is not possible or at least not appropriate from a legal perspective. Furthermore, the clean theory about division of responsibilities between expert and legal representative does not always work this way: there are situations in which the prior actually enters the likelihood ratio, which complicates matters significantly, see e.g.\ \cite{sloomee11}, equations (4.8) and (5.27). Finally, it is highly questionable whether all available evidence is amenable to numerical manipulation.

The theory as developed in this article suggests a slightly different procedure which does not suffer from some of the problems mentioned above and in the introduction.
In our setup, the legal representative still decides on the prior information, be it informative or not. Given this prior, the forensic expert can determine the joint belief structure of $E$ and $G$, and this joint belief structure contains the posterior belief in $G$. Determining the posterior belief, therefore, still is a joint effort of expert and legal representative, but the way in which they execute their roles, is slightly different from the classical case.  

We can illustrate this with the cold case in the island problem discussed in Section \ref{cold}. In the classical case, the expert witness can only deliver the likelihood ratio, which in this case would be $1/p$. Indeed, in the case that the suspect the donor, the probability to find the evidence is 1, and in case he is not, the probability to find this particular profile is nothing but the frequency $p$ of that profile in the population. With uniform prior this leads to the posterior in (\ref{ip2}), but note that other priors would lead to different result. Note that in any case, a prior must be made explicit in order to arrive at posterior probabilities. 

In our setting, if the legal representative confirms that there is no prior information, the expert could simply report (\ref{ip1}) which in that case contains all the available information in the case. If we have no prior information, then no unfounded or subjective choice for a prior needs to be made, apart from the choice of the relevant population. We do not need a uniform prior. If there is prior information, then the expert can set up the corresponding prior belief function, and compute the posterior according to the theory explained above.

We have seen in the examples in the island problem that the numbers obtained are certainly less impressive than in the classical case. This is of course not surprising: starting out with an uninformative prior instead of a uniform prior gives us \emph{less} information to start with. Is this a weakness of the theory of belief functions as we have set out in this paper? We do not think so. It is, in our opinion, better to have a less impressive number which is well founded and not so easy to challenge, than to have a more impressive number which may depend on unfounded arguments or assumptions, and which is easy to dismiss by, say, the defense. As such, our theory may even help to convict actual criminals, even though we insist that if there is no prior information about a certain quantity, we should not pretend there is.

In this article we have described the basic principles of our new theory, with basic forensic examples. Together with the more mathematical companion paper \cite{kerkmee15}, this gives a solid foundation of the theory. Obviously there is a lot of work to be done in the theoretical development and in applications to forensic examples.

\medskip\noindent
{\em Acknowledgment:} We thank Klaas Slooten and Marjan Sjerps for interesting and useful discussions and valuable remarks about earlier versions of this manuscript.

\end{document}